\documentstyle[12pt]{article}
\newcommand{\beq}{\begin{equation}}
\newcommand{\eeq}{\end{equation}}
\newcommand{\beeq}{\begin{eqnarray}}
\newcommand{\eeeq}{\end{eqnarray}}
\newcommand{\eins}{1\hspace{-3pt}{\rm I}}
\begin{document}
\begin{titlepage}
\thispagestyle{empty}

\begin{flushright}
\vspace*{.5cm}
{\large SPBU-IP-97-30
\\ q-alg 9801xxx}
\end{flushright}
\vspace*{.5cm}

\begin{center}
{\Large \bf On Gauss decomposition of quantum groups}\\[.5cm]
{\Large \bf and Jimbo homomorphism{}\footnote
{This research was supported by RFFI
grants No 97-01-01152.}}

\vspace{1cm}

{\large \bf M. A. Sokolov
\footnote{E-mail adress: sokolov@pmash.spb.su}}

\vspace{1cm}

St Petersburg Institute of Mashine Building, \\
Poliustrovskii pr 14, \\
195108, St Petersburg, Russia \\[.5cm]

\vspace{2cm}

\begin{abstract}
It is shown that the properties of the Gauss decomposition of quantum
groups and the known Jimbo homomorphism permit us to realize these
groups as subalgebras of well defined algebras constructed from
generators of the corresponding undeformed Lie algebras.
\end{abstract}
\end{center}
\end{titlepage}

\newpage

{\bf 1.} The purpose of this Letter is to show that every quantum group
\cite{D,J,FRT1} from the Cartan list can be considered as a subalgebra
of a tensor product of well defined algebras constructed from the
corresponding classical Lie algebra.  Such a consideration is based on
the quantum algebra homomorphism constructed by Jimbo in his ground
work \cite{J} and specific properties of the Gauss decomposition of
the standard quantized groups of the series $A_n,\,B_n,\, C_n,\,D_n$
\cite{DS,DKS}.  Recall some relevant definitions and results.  Below
we use the $R-$matrix approach to the theory of quantum groups and
algebras \cite{FRT1}.

According to the $FRT-$ approach the matrix equation

\beq\label{1.1}RT_1T_2=T_2T_1R \eeq
produces the homogeneous quadratic relations for $n^2$ generators
$T=(t_{ij}),\,i,j=1,...,n$ of a unital associative algebra, a number
$n^2$ by $n^2$ matrix $R$ satisfies to the Yang-Baxter equation,
$T_1=T \otimes \eins, \quad T_2= \eins \otimes T$, where $\eins$ is
a unit matrix. Really a number of independent generators is less than
$n^2$ because of the quantum determinant condition for the quantum
groups of the series $A_n$

\beq\label{1.2}
det_qT = \sum_{\sigma}(-q)^{l(\sigma)} \prod_{i=1}^{n} T_{i\sigma (i)}
= 1
\eeq
where the sum is over all the permutations $\sigma$ of the set
$(1,2,...,n)$, and $\sigma(l)$ is length of $\sigma$, and the
additional condition for the orthogonal and symplectic quantum groups

\beq\label{1.3}TCT^tC^{-1}=CT^tC^{-1}T=\eins\eeq
where $T^t$ is the matrix transposed to $T$, and $C$ is a fixed number
matrix \cite{FRT1}.

The quadratic relations for generators of the quantum algebra, dual to
some quantum group, has the form
\beq\label{1.4}R^{(+)}L_{1}^{(\pm)}L_{2}^{(\pm)}=
L_{2}^{(\pm)}L_{1}^{(\pm)}R^{(+)}\eeq

\beq\label{1.5}R^{(+)}L_{1}^{(+)}L_{2}^{(-)}=
L_{2}^{(-)}L_{1}^{(+)}R^{(+)}\eeq
where $L^{\pm}$ are upper- and lower- triangular matrices of the
quantum algebra generators, $R^{(+)}=PRP$, $P$ is the transposition
operator with the properties $P^2=\eins$, $ PA_1P=A_2$, for any matrix
$A$.

The Gauss decomposition for a quantum group is defined as a
transition from the base $T=(t_{ij})$ of the original generators to
the new base \cite{DS}, \cite{DKS}

\beq \label{1.6}T=T_LT_DT_U\eeq
where $T_L=(l_{ik})$ and $T_U=(u_{ik})$ are respectively strictly
lower- and upper- triangular matrices, and $T_D={\rm diag}(d_{kk})$ is
a diagonal matrix. Instead $(T_L,T_D,T_U)$ it is often useful to deal
with the base $(T^{(-)},T^{(+)})$, where $T^{(-)}=T_LT_D,\,T^{(+)}=
T_DT_U$. In the last case we identify the diagonal elements
$t_{ii}^{(-)}=t_{ii}^{(+)}$. Below we shall call the generators of
the both bases by Gauss generators. The commutation relations for these
generators are quadratic too

\beq\label{1.7}RT_1^{(\pm)}T_2^{(\pm)}=T_2^{(\pm)}T_1^{(\pm)}R\eeq
\beq\label{1.8}R_dT_1^{(+)}T_2^{(-)}=T_2^{(-)}T_1^{(+)}R_d\eeq
\beq\label{1.9}R_dT_{D1}T_2^{(-)}=T_2^{(-)}T_{D1}R_d\eeq
\beq\label{1.10}R_dT_1^{(+)}T_{D2}=T_{D2}T_1^{(+)}R_d\eeq
\beq\label{1.11}T_{L1}T_{U2}=T_{U2}T_{L1}\eeq
where $R_d$ is the diagonal part of the $R-$matrix. The complete set
of the equations for Gauss generators see in \cite{DKS}. Here we
remark
only that the requirement (\ref{1.2}) for the quantum special linear
groups takes the form $\prod_{i=1}^n t_{ii}^{(\pm)}\,=\,1$ and the
condition (\ref{1.3}) for orthogonal and symplectic groups is
satisfied by $T^{(\pm)}$.
It is
worth noting that the number of independent Gauss generators for each
orthogonal or symplectic quantum group exactly corresponds to the
number of generators of related classical Lie group.  This is because
the additional deformed relation entering in the quantum group (and
algebra) definition (\ref{1.3}) can be resolved explicitly in the
Gauss base.  Note, as well, that a triangulation procedure \cite{DKS}
can leads to a Gauss base with nonquadratic relations (for example,
this is true for Jordanian quantum group $GL_h(2)$ \cite{DMMZ}).

When the main minors of a $T-$matrix are invertible (in particular,
the element $t_{11}$ is invertible) the set of the Gauss generators is
equivalent, in algebraic structure, to the set of the original ones.
Moreover, in this case we can, in principle, extend the standard
comultiplication operation to the Gauss generators by the homomorphism
property. Below, however, the problems connected with the Hopf algebra
structure of quantum qroups with the Gauss base are not considered.

In the work \cite{FRT2} the following form decomposition of quantum
group generators were introduced

\beq\label{1.12}T=M^{(-)}M^{(+)}\eeq
\beq\label{1.13}[M_1^{(-)},M_2^{(+)}]=0.\eeq
If $M^{(-)}, M^{(+)}$ fulfill the $RTT-$equations (\ref{1.1})
separately then $T=M^{(-)}M^{(+)}$ fulfills as well and, under evident
suggestions \cite{DKS}, conversely. Note that above defined
$T^{(\pm)}$ are not isomorphic to FRT's $M^{(\pm)}$, but we can to
achieve an isomorphism by setting formally
$M^{(-)}=T^{(-)}(\otimes)\eins$, $M^{(+)}=\eins(\otimes) T^{(+)}$. The
symbol $(\otimes)$ means usual matrix multiplication with tensor
product of matrix elements. Thus, we get an isomorphism
$M^{(-)}M^{(+)}\,\simeq\, T^{(-)}(\otimes)T^{(+)}$.

{\bf 2.} The quadratic relations for generators of a quantum algebra
(\ref{1.4}) can be rewritten in the form

$$ R\,\big(L_{1}^{(\pm)}\big)^{-1}\,\big(L_{2}^{(\pm)}\big)^{-1}=
\big(L_{2}^{(\pm)}\big)^{-1} \,\big(L_{1}^{(\pm)}\big)^{-1}\,R.$$
Here we take into account the triangularity of the $q$-matrices
$L^{(\pm)}$ and their invertibility. The last is true because in the
FRT's definition of quantum algebras \cite{FRT1} there exists the
additional requirement for diagonal elements

\beq\label{2.1} {\rm diag}L^{(+)}{\rm diag}L^{(-)}=\eins.\eeq
For that reason any element of $(L^{(\pm)})^{(-1)}$ can be written as
a polynomial in algebra generators $l^{(\pm)}_{i\,j}$ , $i\neq j$ and
inverse diagonal generators $\left(l^{(\pm)}_{i\,i}\right)^{-1}$.

As the additional relations (\ref{1.3}) has the same form for both
orthogonal and symplectic algebras $L^{(\pm)}$  and quantum groups
in the Gauss base $T^{(\pm)}$, we get an algebra (but not a Hopf
algebra) isomorphism

$$T^{(\pm)}\, \simeq \,\big(L^{(\pm)}\big)^{-1}.$$
Under this isomorphism to the elements
$t_{i,\,i+1}^{(+)}$~, $t_{i+1,\,i}^{(-)}$ of Borel matrices
$T^{(\pm)}$  correspond the elements $l_{i,\,i+1}^{(+)}$~,
$l_{i+1,\,i}^{(-)}$   of the matrices $L^{(\pm)}$:

$$ t_{i,\,i+1}^{(+)} \leftrightarrow -\,\big(l_{i\,i}^{(+)}\big)^{-1}\,
l_{i,\,i+1}^{(+)}\, \big(l_{i+1,\,i+1}^{(+)}\big)^{-1}; $$

$$ t_{i+1,\,i}^{(-)} \leftrightarrow -\, 
\big( l_{i+1,\,i+1}^{(-)}\big)^{-1}\,
l_{i+1,\,i}^{(-)}\,\big( l_{i\,i}^{(-)}\big)^{-1} . $$
(except the algebras of the series $D_n$, see \cite{FRT1,DS2}). The
elements $l_{i,\,i+1}^{(+)}$~, $l_{i+1,\,i}^{(-)}$ are associated with
the simple roots of classical algebras through known identification
{\cite{FRT1}. As a consequence we obtain that  the elements
$t_{i,i+1}^{(+)},\, t_{i+1,i}^{(-)}$ satisfy the identities
which slightly differ from deformed Serre identities (\cite{D,J}).
These identities together with the commutation
relations between the generators $t_{i,i+1}^{(+)},\,t_{i+1,i}^{(-)}$
and the diagonal elements $t_{ii}^{(+)}, t_{ii}^{(-)}$ produce the set
of formulas which completely defines the quantum group in Gauss
base.

{\bf 3.} Hereafter we restrict our consideration to the quantum groups
$SL_q(n)$, but all the formulas can be easily extended to the
orthogonal and symplectic cases.

Let $(H_i,X_i^{(\pm)}),\,(i=1,....,n-1)$ be the Chevalley base for
the classical Lie algebra $sl(n,C)$. The commutation relations of
the Borel subalgebra $\tt b^{\pm}$ generators and the Serre identities
are defined by the algebra simple positive roots $\alpha_i $:

\beq\label{3.1}[H_i,H_j]=0, \quad [H_i,X_j^{(\pm)}]=\pm
(\alpha_i,\alpha_j)X_j^{(\pm)}, \quad i,j=1,....,n ,\eeq

\beq\label{3.2}(adX_i^{(\pm)})^{1-A_{ij}}(X_j^{(\pm)})=0,  i\not=j,\eeq
where  $A$ is the Cartan matrix
$A_{ij}=2(\alpha_i,\alpha_j)/(\alpha_j,\alpha_j).$

Let us consider $n$ linear dependent elements from the Cartan
subalgebra of $sl(n,C)$ defined as

\beq\label{3.3}\widetilde H_i=n^{-1}\big(\sum_{k=i}^{n-1}(n-k)H_k \,
-\, \sum_{k=1}^{i-1}kH_k \big),\quad i=1,...,n\eeq
Such elements were introduced in \cite{FRT1} to establish the relation
between the Drinfeld-Jimbo and $R-$matrix formulations of the quantum
algebra theory. Note, however, that $\widetilde H_i$ were regarded as
deformed objects in the cited paper.

Introduce the elements $K_i=e^{h\widetilde H_i}$ where $q=e^h$ is a
deformation parameter. In fact, $K_i$ are the elements of the maximal
commutative subgroup of the complex Lie group $SL(n,C)$. Thus, we have
well defined adjoint action of $K_i$ on $sl(n,C)$

\beq\label{3.4}
AdK_i(X_j^{(\pm)})=\left\{\begin{array}{ll}
q^{\pm1}X_j^{(\pm)} & {\rm if} \ j=i \\
q^{\mp1}X_j^{(\pm)} & {\rm if} \ j=i \mp 1 \\
X_j^{(-)} & {\rm if} \ j \not=i,i \mp 1
\end{array}
\right.
\eeq
The last formulas can be easily proved by direct calculations using
(\ref{3.1},\ref{3.3}).

Denote by ${\cal U}^{+}_i$ an unital algebra generated by
$\big(K_i^{\pm1},\, X_i^{(+)}\big)$ and by ${\cal U}^{-}_i$
an unital algebra generated by $\big(K_i^{\pm1},\, X_i^{(-)}\big),\;$
$${\cal U}^{\pm} \,=\,{\cal U}_1^{\pm} \otimes {\cal U}_2^{\pm}
\otimes...\otimes {\cal U}_{n-1}^{\pm}.$$ In view of (\ref{3.4})
${\cal U}^{\pm}_i$ is a subalgebra of $G_{\cal H} \otimes
{\cal U}^{\pm}(X^{\pm}_i)$, where $G_{\cal H}$ is the maximal
commutative subgroup of $SL(n,C)$ and ${\cal U}^{\pm}(X^{\pm}_i)$
is the algebra generated by a single element. Consider,
following Jimbo \cite{J}, the homomorphism $\delta^{(\pm)}$ :
${\cal U}^{\pm}_i\,\longrightarrow\,{\cal U}^{\pm}$

\beq\label{3.5} \delta^{(\pm)}(K_i)\,=\,t_{ii}^{(\pm)}
\,=\, \big(K_i^{\pm1}\big)^{\otimes(n-1)} \quad i=1,...,n\eeq

\beq\label{3.6} \delta^{(+)}(X_i^{(+)})\,=\,t_{i,i+1}^{(+)}=
f_i K_i^ {\otimes (i-1)}\otimes X_i^{(+)} 
\otimes K_{i+1}^{\otimes (n-i-1)}
\quad i=1,....,n-1\eeq

\beq\label{3.7} \delta^{(-)}(X_i^{(-)})\,=\,t_{i+1,i}^{(-)}\,=\,
g_i \;\big(K_i^{-1}\big)^{\otimes(i-1)} \; \otimes \; X_i^{(-)} \;
\otimes\big(K_{k+1}^{-1}\big)^{\otimes (n-i-1)}\quad i=1,....,n-1.\eeq
In the above formulas $K^{\otimes n}\,=\, K \otimes K \otimes \ldots
\otimes K,$ ($n$ tensor multipliers);\,   $K^{\otimes 1}=K$;\,
$f_i,g_i$ are arbitrary constants. The homomorphism $\delta^{(\pm)}$
slightly differs from the Jimbo's one \cite{J} in its form but we
stress here again that it is a homomorphism of undeformed objects.

Consider the elements $t_{i,i+1}^{(+)}$ defined by (\ref{3.5}).
In view of (\ref{3.4}), they satisfy the following relations
$$t_{i,i+1}^{(+)}t_{j,j+1}^{(+)}=\left\{ \begin{array}{ll}
q^{(\pm 2)} t_{j,j+1}^{(+)}t_{i,i+1}^{(+)} & if \ j=i \pm 1 \\
t_{j,j+1}^{(+)}t_{i,i+1}^{(+)} & if \ j \not= i \pm 1 \end{array}
\right. $$
Put $X_i^{(+)}\,=\,t_{i,i+1}^{(+)}$. Then the deformed Serre
identities, which have the form in our case
$$(X_i^{(+)})^2X_j^{(+)})-q^{\pm1}(q+1/q)X_i^{(+)})X_j^{(+)})
X_i^{(+)})+q^{\pm 2}X_j^{(+)})(X_i^{(+)}))^2=0$$
for $j=i \pm 1,$ and the commutation relations
$$[X_i^{(+)},X_j^{(+)}]=0 \quad for \, j \not= i \pm 1$$
are satisfied by these $t_{i,i+1}^{(+)}$.
Therefore, one can uniquely reconstruct the other elements of
the $T^{(+)}$ by the formula

\beq\label{3.8} t_{i,i+k}^{(+)}=\lambda^{1-k}\big(\prod_{l=1}^{k-1}
(t_{i+l,i+l}^{(+)})^{-1}\big)[t_{i,i+1}^{(+)},[t_{i+1,i+2}^{(+)},[.....
[t_{i+k-2,i+k-1}^{(+)},t_{i+k-1,i+k}^{(+)}].....]  \eeq
which easily can be verified by induction.
Application of (\ref{3.8}) to  (\ref{3.5},\ref{3.6}) yields

\beq\label{3.9} t_{i,k+1}^{(+)}=\prod_{j=i}^k f_j \; K_i^{\otimes
(i-1)} \; \otimes \;\big(\bigotimes_{l=i}^k X_l^{(+)}\big)\; \otimes
K_{k+1}^{\otimes (n-i-1)}\eeq
where $\bigotimes_{l=i}^k F_l=F_i \otimes F_{i+1}\otimes....\otimes
F_k$,\, $\bigotimes_{l=i}^i F_l=F_i$,\,$ k \geq i$.
It easy to verify that the constructed elements
$T^{(+)}=(t_{ij}^{(+)}),\, i \leq j$ satisfy $RTT-$equation
(\ref{1.1}).

The reconstruction of the elements of $T^{(-)}$ can be carried out in
close analogy with the above $T^{(+)}$- case. Note, however, that
$\delta^{(+)}=\kappa \circ \delta^{(-)}$ where  $\kappa$ is
the Cartan automorphism
$$ H_i \rightarrow - H_i; \quad X_l^{(+)} \rightarrow X_l^{(-)}.$$
As a result we have immediately

\beq\label{3.10}t_{k+1,i}^{(-)}=\prod_{j=i}^k g_j
\;\big(K_i^{-1}\big)^{\otimes (i-1)} \;\big(\bigotimes_{l=i}^k
X_l^{(-)}\big) \; \otimes \big(K_{k+1})^{-1}\big)^{\otimes (n-i-1)}.
\eeq
The above realization of $T^{(\pm)}$ elements permits us accomplish
the final step which consists in constructing of the matrix
$T\,=\,T^{(-)}(\otimes)T^{(+)}$. The latter means that the generators
of the quantum group $SL_q(n)$ and all its elements belong to the
algebra:

\beq\label{3.11} {\cal U}\,=\,{\cal U}^- \otimes {\cal U}^+\,=\,
\bigotimes_{i=1}^{n-1}{\cal U}_i^{-}
\otimes
\bigotimes_{j=1}^{n-1}{\cal U}_j^{+}\eeq

{\bf 4.} To illustrate the above construction let us consider the
quantum group $SL_q(2)$ as an example. The Gauss decomposition for
this group is of the form

\begin{eqnarray*}\label{4.1}
T &=&\left(
\begin{array}{cc}
t_{11} & t_{12} \\
t_{21} & t_{22}
\end{array}
\right) = \left(
\begin{array}{cc}
1 & 0 \\
l & 1
\end{array}
\right)  \left(
\begin{array}{cc}
A & 0 \\
0 & A^{-1}
\end{array}
\right)  \left(
\begin{array}{cc}
1 & u \\
0 & 1
\end{array}
\right) = \left(
\begin{array}{cc}
A  & Au \\
lA & lAu+A^{-1}
\end{array}
\right).
\end{eqnarray*}
Commutation relations between the Gauss generators $A,A^{-1},l,u$
are very simple
\beq\label{4.2}
Au=quA,\,Al=qlA,\,uA^{-1}=qA^{-1}u,\,lA^{-1}=qA^{-1}l,\,[u,l]=0.
\eeq
According to the above prescription, introduce the matrices
$T^{(\pm)}:$

\begin{eqnarray*}\label{4.3}
T^{(-)}\;=\;\left(
\begin{array}{cc}
A & 0 \\
lA & A^{-1}
\end{array}
\right)\,=\,\left(\begin{array}{cc}
q^{-H/2} & 0 \\
gX^{(-)} & q^{H/2}
\end{array}\right);
\end{eqnarray*}

\begin{eqnarray*}\label{4.4}
T^{(+)}\; =\; \left(
\begin{array}{cc}
A & Au \\
0 & A^{-1}
\end{array}
\right)\,=\,\left(\begin{array}{cc}
q^{H/2} & fX^{(+)}\\
0 & q^{-H/2}\end{array}\right).
\end{eqnarray*}

Finally, we have

\begin{eqnarray*}\label{4.5}
T\; =\;T^{(-)} \otimes T^{(+)} \, = \, \left(\begin{array}{cc}
q^{-H/2} \otimes q^{H/2} & fq^{-H/2}\otimes X^{(+)}\\
gX^{(-)} \otimes q^{H/2}& q^{H/2}\otimes q^{-H/2}+fgX^{(-)}\otimes
X^{(+)} \end{array}\right).
\end{eqnarray*}

{\bf 5.} Let us make some remarks in conclusion.
\begin {itemize}
\item The above construction can be easily generalized to a
multiparameter case. For example, for the quantum group
$GL_{p,q}(2)$, defined by the $R$-matrix 

\begin{eqnarray*}
R_{p,q}=\left(
\begin{array}{cccc}
k & 0 & 0 & 0\\
0 & p & 0 & 0\\
0 & k-pq/k & q & 0\\
0 & 0 & 0 & k
\end{array}
\right),
\end{eqnarray*}
we have the following construction
$$T=\left(
\begin{array}{cc}
(k/p)^{-H/2} & 0 \\
c_{-}X^{(-)} & (k/q)^{H/2}
\end{array}\right) (\otimes) \left(
\begin{array}{cc}
(k/q)^{H/2}& c_{+}X^{(+)} \\
      0    & (k/p)^{-H/2}
\end{array}
\right)$$
where $c_{\pm}$ are arbitrary constants.
The one-parameter case corresponds to the change
$k \rightarrow q^{1/2} , \quad p \rightarrow q^{-1/2}
\quad q \rightarrow q^{-1/2}.$
\item Using any matrix representation of the algebra $sl(n,C)$, one
can obtain a matrix representation for the corresponding quantum group.
For instance, using the lowest dimensional representation of $sl(2)$
by the matrices
$X^{(+)}=\left(
\begin{array}{cc}
0&1\\0&0
\end{array}\right); \quad
X^{(-)}=\left(
\begin{array}{cc}
0&0\\1&0
\end{array}\right); \quad
H=\left(
\begin{array}{cc}
1&0\\0&-1
\end{array}\right),$
one has for $SL_q(2)$ in view of (\ref{4.5})

$$t_{11}=\left(
\begin{array}{cccc}
1 & 0 & 0 & 0\\
0 & q & 0 & 0\\
0 & 0 & q^{-1} & 0\\
0 & 0 & 0 & 1
\end{array}\right); \quad
t_{12}=\left(
\begin{array}{cccc}
0 & 0 & 1 & 0\\
0 & 0 & 0 & q\\
0 & 0 & 0 & 0\\
0 & 0 & 0 & 0
\end{array}
\right); \quad
t_{21}=\left(
\begin{array}{cccc}
0 & 0 & 0 & 0\\
1 & 0 & 0 & 0\\
0 & 0 & 0 & 0\\
0 & 0 & q^{-1} & 0
\end{array}
\right);$$
$$t_{22}=\left(
\begin{array}{cccc}
1 & 0 & 0 & 0\\
0 & q^{-1} & 1 & 0\\
0 & 0 & q & 0\\
0 & 0 & 0 & 1
\end{array}
\right).$$
\item Generalizing the matrix representation of the above subsection
we can formally associate with any representation of $sl(n,C)$
a representation of the corresponding quantum group $SL_q(n)$.
This construction is under investigation.
\item We can decrease a number of tensor multipliers in the considered
construction if we use the dual algebra $sl^{*}(n)$ rather than
$sl(n)$. In the case of $sl^{*}(2)$ we have the commutation relations
for its generators $\widetilde H,\widetilde X^{\pm}:$

\beq\label{5.1} [\widetilde H, \widetilde X^{\pm}]= \widetilde
X^{\pm},\, [\widetilde X^+,\widetilde X^-]=0.\eeq

The last commutator suggests to use the Gauss decomposition in the
form (\ref{1.6}). As a result we have the realization

$$T =\left(
\begin{array}{cc}
1 & 0 \\
c_- \widetilde X^- & 1
\end{array} \right)
\left(
\begin{array}{cc}
q^{\widetilde H} & 0 \\
0 & q^{-\widetilde H}
\end{array}\right) \left(
\begin{array}{cc}
1 & c_- \widetilde X^{+}\\
0 & 1
\end{array}\right)$$
where $c_{\pm}$ are arbitrary constants.
\item From the explicit form of the $SL_q(n)$ generators we can
formally obtain a classical limit of $q$-matrix $T$
by direct
differentiation of $t_{ij}$ in $h$ and setting $h=0$.
For $SL_q(2)$,
choosing the parameters $f=q^{-1}\lambda ,\, g=-q \lambda$
\cite{DS2} we get

$$M\,=\,dT/dh(h=0)\; =\;
\left(\begin{array}{cc}
1/2(1 \otimes H - H \otimes 1) &  2 (1 \otimes X^{(+)})\\
-2(X^{(-)} \otimes 1)& -1/2(1 \otimes H - H \otimes 1)
\end{array}\right).$$
The elements of $M$-matrix satisfy the commutation relations
(\ref{5.1}).
\end{itemize}

{\bf Acknowledgments.} The author would like to thank E.V.
Damaskinsky, P.P. Kulish and V.D.  Lyakhovsky for useful discussions.
The work was
supported in part by the Russian Fund for Fundamental Research (grant
No 97-01-01152).

\end{document}